\documentstyle[amscd,12pt,psamsfonts]{amsart}

\newtheorem{defn}{Definition}
\newtheorem{thm}[defn]{Theorem}

\newtheorem{lem}[defn]{Lemma}

\theoremstyle{remark}
\newtheorem{rem}{Remark}
\theoremstyle{remark}

\numberwithin{equation}{section}
\numberwithin{defn}{section}

\begin{document}


\newcommand\spk{{\operatorname{Spec}}(k)}
\renewcommand\sp{\operatorname{Spec}}       
\renewcommand\sf{\operatorname{Spf}}       
\newcommand\proj{\operatorname{Proj}}
\newcommand\aut{\operatorname{Aut}}
\newcommand\grv{{\operatorname{Gr}}(V)}
\newcommand\gr{\operatorname{Gr}}
\newcommand\glv{{\operatorname{Gl}}(V)}
\newcommand\glve{{\widetilde{\operatorname{Gl}}(V)}}
\newcommand\gl{\operatorname{Gl}}
\renewcommand\hom{\operatorname{Hom}}
\newcommand\End{\operatorname{End}}
\renewcommand\det{\operatorname{Det}}    
\newcommand\tr{\operatorname{Tr}}      
\newcommand\detd{\operatorname{Det}^\ast}
\newcommand\im{\operatorname{Im}}
\newcommand\res{\operatorname{Res}}
\renewcommand\deg{\operatorname{deg}}
\newcommand\id{\operatorname{I}}
\newcommand\rank{\operatorname{rank}}
\newcommand\limi{\varinjlim}
\newcommand\limil[1]{\underset{#1}\varinjlim\,}
\newcommand\limp{\varprojlim}
\newcommand\limpl[1]{\underset{#1}\varprojlim\,}
\newcommand\dk{{\bf\delta}}       

\renewcommand\o{{\mathcal O}}      
\renewcommand\L{{\mathcal L}}       
\newcommand\B{{\mathcal B}}       
\renewcommand\c{{\widehat C}}
\renewcommand\P{{\mathcal P}}    
\newcommand\Q{{\mathcal Q}}    
\newcommand\D{{\mathcal D}}    
\newcommand\C{{\mathbb C}}    
\newcommand\Z{{\mathbb Z}}    
\newcommand\A{{\mathbb A}}    
\newcommand\M{{\mathcal M}}
\newcommand\T{{\mathcal T}}
\newcommand\U{{\mathcal U}}

\newcommand\w{\widehat}
\renewcommand\tilde{\widetilde}  
\newcommand\f{\overset{\to}f}

\newcommand\iso{@>{\sim}>>}
\renewcommand\lim{\limpl{A\in\B}}
\newcommand\beq{       
      \setcounter{equation}{\value{defn}}\addtocounter{defn}1
      \begin{equation}}

\renewcommand{\thesubsection}{\thesection.\Alph{subsection}}

\title{Algebraic Solutions of the Multicomponent KP Hierarchy}
\author[F. J. Plaza Mart\'{\i}n]{F. J. Plaza Mart\'{\i}n\\  \medskip\tiny
Departamento de Matem\'aticas \\ Universidad de Salamanca}

\address{{\it Current address}:\newline DPMMS, University of Cambridge,
Mill Lane 16, Cambridge CB2 1SB, United Kingdom
 \newline \indent {\it Permanent address}:\newline Departamento de
Matem\'aticas, Universidad de Salamanca,  Plaza de la Merced 1-4\\
Salamanca 37008. Spain.}

\thanks{
This work is partially supported by the CICYT research contract n.
PB96-1305 and Castilla y Le\'on regional government contract SA27/98.}

\email{\newline fjpm@@dpmms.cam.ac.uk \newline fplaza@@gugu.usal.es}

\begin{abstract}
It is shown that it is possible to write down tau functions for the
$n$-component KP hierarchy ($n$-KP) in terms of non-abelian theta
functions. This is a generalization of the rank 1 situation; that is,  the
relation of theta functions of Jacobians and tau functions for the KP
hierarchy. 
\end{abstract}

\date{July 12th, 1999}

\maketitle


\setcounter{tocdepth}1
\tableofcontents

\section{Introduction}

This paper is concerned with the generalization for higher rank of the
relation of theta functions of Jacobians and the KP hierarchy. More
precisely, the well known relation between the theta function of a
Jacobian and the tau function of a certain point of the infinite
Grassmannian $\gr^0(k((z)))$ (\cite{Kr2,SW,Sh}, see also
\cite{MP,Pl,PW}) is generalized for non-abelian theta functions and
$\gr^0(k((z))^{\oplus n^2})$. From this
point of view, the main result is the Theorem \ref{thm:tau-theta} that
shows that it is possible to write down tau functions for the $n$-component
KP hierarchy ($n$-KP) in terms of non-abelian thetas. However, further
research must be made to obtain  explicit expressions.

It is worth remarking that this result is twofold; on the one hand, it
shows that non-abelian thetas give solutions for $n$-KP; on the other, it
suggests the possibility of characterize non-abelian thetas
in terms of differential equations of the $n$-KP (see \cite{Sh,MP} for the
rank 1 case). 

In a certain sense, \S2 is the core of the paper. It contains three new
results; namely, Lemma \ref{lem:main} (Addition formula) and
Theorems \ref{thm:BAmain} and
\ref{thm:bilid} (Residue Bilinear Identity). These results are the key for
proving that  there exists a 1-1 correspondence between wave functions for
the
$n$-KP and Baker-Akhiezer (BA) functions of points of the infinite
Grassmannian
$\gr^0(k((z))^{\oplus n^2})$ (see
\S3). Then, our strategy is rather simple (see \S5): a generalization of
the well-known Krichever map (\cite{Kr2,Mum,Mul}) is used to show that
there is a subscheme:
$$\iota:\w\U(r,d)\hookrightarrow\gr^0(k((z))^{\oplus n^2})$$
together with a projection onto the moduli space of vector bundles,
$\pi:\w\U(r,d)\to\U(r,d)$ (some data must be previously fixed) such that:
$$\begin{gathered}
\iota_{F}^*\Omega_+\,=\,\pi^*\theta_F \\
\iota_{F}^*\det\,\iso\,\pi^*\o(\Theta_{[F]})
\end{gathered}$$
(see \S5 for notations and precise statements). This proves that
non-abelian thetas give rise to tau functions for the $n$-KP.

Before explaning how the paper is organized, let us point out an
intermediate result that deserves special mention; namely, Theorem
\ref{thm:equ-mod} that computes the equations defining the subscheme of
$\gr(k((z))^{\oplus n})\times
\gr(k((z))^{\oplus n\cdot r})$ whose set of rational points is:
$$\big\{(C,p_1,\ldots,p_n,\alpha_1,\ldots,\alpha_n,M,\beta)\big\}$$
where $C$ is a curve, $p_i\in C$, $\alpha_i:\w\o_{C,p_i}\simeq k[[z]]$,
$M$ is a rank $r$ torsion free sheaf, and $\beta:\w M_{\{p_i\}}\iso
k[[z]]^{\oplus n\cdot r}$.

In \S2 the approach of \cite{MP} to the KP hierarchy, which is based in
the ``geometry of formal curves'', is developed for the $n$-KP; that
is, for the infinite Grassmannian of $E((z))$ where $E$ is a finite
dimensional $k$-vector space. This enables us to define tau functions for
the $n$-KP in terms of global sections of the determinant bundle and to
generalize the Addition formula for this situation (see Lemma
\ref{lem:main}). Similarly to the rank 1 case, the BA function
of an arbitrary point of $\gr(E((z)))$ is defined as a certain deformation of the tau
function. However, from the Theorem \ref{thm:BAmain} it follows that our
definition agrees with the standard one (\cite{Kr1,Kr2,KrN}, see
\cite{Kr4,Pr,PW} for overviews on the subject) for those  points coming
from algebro-geometric data. This section finishes with the generalization
of the Residue Bilinear Identity (Theorem
\ref{thm:bilid}).

Section \S3 recalls the definition of the $n$-KP and shows that there
exists a 1-1 correspondence between wave functions for the $n$-KP and BA
functions of points of the infinite Grassmannian $\gr^0(k((z))^{\oplus
n^2})$. 

Moduli spaces are studied in section \S4 in terms the Krichever map and
infinite Grassmannians. Although it is not needed for our main purpose,
we have considered convenient to include here the equations for these
moduli spaces (Theorem \ref{thm:equ-mod}).

The last section unveils the deep relation among tau functions of
algebro-geometric points and non-abelian thetas, which is the
``expected'' generalization of the rank 1 situation.

Finally, I hope that these results will help in the study of some related
problems; particularly, those related with higher rank vector bundles over
curves and infinite Grassmannians and with differential operators (e.g.
 higher rank commutative subrings of
differential operators, Darboux-B\"acklund transform, etc).

I would like to express my gratitude to Prof. G. Segal for inviting me to
the DPMMS at University of Cambridge (UK) where this work has been done.

\section{Infinite Grassmannians}
\subsection{Background}

Let us address the reader to \cite{AMP,MP} for the scheme-theoretic
approach to infinite Grassmannians. However, it is
convenient to recall some basic facts in order to fix notations and to
point out the statements which we shall need.

Since we are concerned with the multicomponent KP hierarchy, we shall not
deal with general infinite Grassmannians, but only with that of $V=E((z))$
where $E$ is a $n$-dimensional $k$-vector space. Denote $E[[z]]$ by
$V^+$ and consider the linear topology in $V$ given by $\{z^m V^+\vert
m\in\Z\}$ as a basis of neighborhoods of $(0)$.

Then, we know  that there exists a $k$-scheme $\grv$ locally covered by
the open subschemes:
$$F_A(S)\,:=\, \left\{ \text{sub-$\o_S$-modules $\L\subset\w V_S$ such
that }\L\oplus\w A_S=\w V_S\right\}$$ where $S$ is a $k$-scheme and:
\begin{itemize}
\item $A\subset V$ is a subspace such that $A\sim V^+$; that is,
$\dim(A+V^+)/A\cap V^+<\infty$;
\item $\w A_S$ is defined by $\limpl{B\sim V^+}((A/A\cap B)\otimes\o_S)$
for a subspace $A\subseteq V$. (For instance, when $E=k$ one has $\w V_S^+=
\limpl{m}\o_S[z]/z^m=:\o_S[[z]]$, and $\w V_S=\limil{m}z^{-m}\o_S[[z]]=:
\o_S((z))$).
\end{itemize}

The generalization of some ``good'' properties of the $1$-dimensional case
requires a choice of a chain of strict inclusions $V^+_0=V^+\subset
V^+_1\subset\ldots\subset V^+_n =z\cdot V^+$ (note that it follows that
the inclusions are of codimension
$1$). Some remarkable facts are:
\begin{itemize}
\item the function:
$$\begin{aligned}
\grv &\longrightarrow \Z \\ L &\mapsto \dim(L\cap V^+)-\dim(V/L+V^+)
\end{aligned}$$ gives the decomposition of $\grv$ in connected components,
which will be denoted by $\gr^n(V)$ ($n\in\Z$);
\item there is a natural line bundle on $\grv$ defined (on the connected component $\gr^n(V)$) by  the determinant of the perfect complex:
$$\L\oplus (\w V^+_m)_{\o_{\gr^m(V)}}\longrightarrow (\w
V_m)_{\o_{\gr^m(V)}}$$ where
$\L$ is the submodule of $\w V_{\o_{\gr^m(V)}}$ corresponding to the
universal object, and $V_m$ is defined as $z^q\cdot V^+_r$ with $m=q\cdot
n+r$ and $0\leq r<n$;
\item the addition morphism in the above complex gives canonically a
global section of the dual of that bundle, $\Omega_+\in H^0(\grv,\detd)$.
\end{itemize}

\begin{rem} Although as abstract scheme $\grv$ is independent of the
dimension of $E$, it is straightforward that the groups acting (naturally)
on the Grassmannians do depend on it. The standard procedure to introduce
these Grassmannians consists of a ``re--labelling'' the indexes of
$k((z))$ and taking $V^+_m =T^m(V^+)$, where:
$$T\,:=\,\pmatrix 0 & 1 & 0 & \ldots &0 \\
0 & 0 & 1 &\ldots &0 \\
\vdots & \vdots & &\ddots & \vdots \\
0 & 0 & \ldots  & 0 & 1 \\
z & 0 & \ldots & 0 & 0
\endpmatrix$$
but in this case the induced $k((z))$-module structures are different.
\end{rem}

For simplicity's sake, we will assume that:
\begin{itemize}
\item  $E=k^{\oplus n}$; that is, we consider a basis $\{e_1,\ldots,e_n\}$ of $E$, so that the elements of $E((z))$ may be thought as $E$-valued series in $z$ or as $n$-tuples of series;
\item  $k$ is an algebraically closed field of characteristic $0$.
\end{itemize}

\subsection{$\tau$ Functions}

Analogously to the one dimensional case, the formal geometry language (\cite{AMP,MP}) will be the base of our approach to the definition
$\tau$-function. From now on, a pair $(E,\T)$ consisting of a $n$-dimensional $k$-vector
space and a semisimple commutative subalgebra $\T\subseteq \End(E)$
will be fixed. 

Although the constructions and results below hold in greater generality
(e.g. $\dim \T\leq \dim E$), we will assume that $\T$ is $n$-dimensional
and that there exist a basis
$\{T_1,\ldots ,T_n\}$ such that
$T_i(e_j)=\delta_{ij}e_i$ for $1\leq i\leq n$ and $1\leq j\leq n$. 
(The existence follows from \cite{Bou} Chp.~8,\S9 n.3). 

Motivated by the fact that the choice of a basis of $\T$ induces an isomorphism of the completion of the symmetric
algebra generated by $\T$, and $k[[t_1,\dots,t_n]]$, define the
$n$-dimensional formal variety by:
$$\w C^n\,:=\, \sf(k[[t_1,\dots,t_n]])$$
And let $\Gamma_-^n$ be the direct limit of the symmetric products of $\w
C^n$, $\limil{r}S^r\w C^n$.

Now, we will study an action of $\Gamma_-^n$
that, roughly said, is induced by that of $\T$ on $V$ at an infinitesimal
level.

Analogously to the case $n=1$ (see Theorem~3.6 of
\cite{AMP}), one has:

\begin{thm}  
$\Gamma_-^n$, is a formal group scheme whose
$S$-valued points are $n$-tuples of series on $z^{-1}$:
{\small $$\big(1+\sum_{i>0} a_{i1}z^{-i}, \ldots,
1+\sum_{i>0} a_{in}z^{-i}\big)$$}
(where $a_{ij}\in H^0(S,\o_S)$ are nilpotents) with componentwise
multiplication as composition law.
\end{thm}

Let $\prod^r\w C^n$ be the formal spectrum of
$ k[[\{t_{i j}\}_{\Sb  1\leq j\leq n\\  1\leq i\leq r\endSb}]]$. 
Then, the ring of $\Gamma_-^n$ is $\limpl{r>0} k[[\{s_{i j}\}_{\Sb 
1\leq
 j\leq n\\ 1\leq i\leq r\endSb} ]]$ (where $s_{i j}$ is to be
understood as the $i$-th symmetric function on
$t_{1 j},t_{2 j},\ldots$). Since $\operatorname{char}(k)=0$, the
exponential map gives an isomorphism of $\Gamma^n_-$ with an additive
group scheme, the universal element might be written as
$(\exp(\sum_{i>0}s_{i1}z^{-i}), \ldots,
\exp(\sum_{i>0} s_{in}z^{-i}))$.

Consider the action $\mu :\Gamma_-^n\times \grv\to\grv$ induced by that 
of $\Gamma_-^n$ on $V$ given by componentwise multiplication; or,
equivalently by:
{\small $$
(\exp(\sum_{i>0}s_{i1}z^{-i}), \ldots, \exp(\sum_{i>0} s_{in}z^{-i}))\cdot
v\,:=\, \exp\big(\sum_{j=1}^n\sum_{i>0} s_{ij}z^{-i}T_j\big)(v)$$}

Note that $\mu $ preserves the determinant bundle; that $(\mu_U^n)^*\detd$
is trivial; and, that the group structure of $\Gamma_-^n$ gives a
trivialization of it.

\begin{defn}\label{defn:tau}
Then $\tau$-function of a rational point $U\in\grv$, $\tau_U$, is the
inverse image of the global section of the determinant bundle
$\Omega_+$ by the morphism $\mu_U:\Gamma_-^n\times\{U\}\to \grv$.
\end{defn}

\subsection{Baker-Akhiezer Functions}

Recall that the Abel morphism: 
$$\sf k[[\bar z]]\times\Gamma_-\to\Gamma_-$$
is the associated with the series $(1-\frac{\bar z}z)^{-1}\cdot
(1+\sum_{i>0}a_iz^{-i})$. Let $\phi_j:\sf(k[[\bar z]])\times
\Gamma_-^n\to \Gamma_-^n$ be the morphism given by the Abel morphism
in the $j$-th entry and by the identity on the others.

\begin{defn}\label{defn:BA}
The Baker-Akhiezer function of a rational point $U\in\gr^r(V)$ (with
$\Omega_+(U)\neq 0$),
$\psi_U(z,s)$, is the vector valued function:
$$\psi_U(z,s)\,:=\,
\exp(-\sum_{i,j}T_j\frac{s_{i j}}{z^i})\cdot
\frac1{\tau_U(s)}\cdot
(\phi_1^*(\tau_U),\dots,\phi_n^*(\tau_U))$$
\end{defn}

\begin{rem}
In our definitions of tau and BA functions, the commutativity of $\T$ is extremely important. It implies that 
expressions like $\prod_{\Sb 1\leq i\leq N\\ 1\leq j\leq n\endSb}
(1-T_j\frac{t_{i j}}z)$,
$\exp(-\sum_{i,j}T_j\frac{s_{i j}}{z^i})$ are well defined
and, further, that the Abel morphisms $\phi_j$ and the morphism
$\mu_U$ are compatible. 
\end{rem}

In order to introduce the adjoint BA function, one must assume that there
is a non-degenerate symmetric pairing:
$$T_2:E\times E\to k$$

Then, $E((z))$ carries a natural non-degenerate pairing, $\res$, given by:
{\small $$\res (\sum f_i z^i)\cdot(\sum g_j z^j)\,:=\,
\res_{z=0}\sum_{i,j}T_2(f_i,g_j)z^{i+j}{dz}
\,=\, \sum_i T_2(f_{i},g_{-i-1})$$}

When $E$ is the algebra of matrices, $M_{r\times s}(k)$ (including the
case $r=1$), we shall consider the pairing given by:
$$(A,B)\,\mapsto
\tr(A\cdot B^t)$$

\begin{defn}
The adjoint BA function of a rational point $U\in\grv$ is:
$$\psi_U^*(z,s)\,:=\, \psi_{U^\perp}(z,-s)$$
\end{defn}

\subsection{Addition formula}

Let $N>0$ be an integer. Let $U\in\grv$ be a rational point. Consider
the morphism:
$$\bar\mu_U^N:\prod^N\w C^n\longrightarrow \grv$$
given by:
$$\big(\prod_{\Sb 1\leq i\leq N\\ 1\leq j\leq n\endSb}
(1-T_j\frac{t_{i j}}z)\big)^{-1}(U)$$ (the ring of the $j$-th
copy of
$\prod^N\w C$ in $\prod^N\w C^n$ is
$k[[t_{1 j},\ldots, t_{N j}]]$).

\begin{lem}\label{lem:main}
Let $U\in\gr^0(V)$ satisfy $\Omega_+(U)\neq 0$. Then, for all $N>>0$ the
inverse image of $\Omega_+$ by $\bar \mu_U^N$ is given (up to a non-zero
scalar) by the following expression:
{\small
$$\big(\prod_{j=1}^n\Delta_j\big)^{-1}\cdot
\operatorname{det}
\pmatrix f^1_1(t_{11}) & \ldots & f^1_1(t_{N1})  & \ldots & f^1_n(t_{1n})
& \ldots & f^1_n(t_{Nn}) \\
\vdots & & & & & & \vdots \\
  f_1^M(t_{11}) & \ldots & f_1^M(t_{N1})  & \ldots & f_n^M(t_{1n})
& \ldots & f^M_n(t_{Nn}) 
\endpmatrix$$}
 where $M:=N\cdot n$ and $\{f^i=(f^i_1,\ldots, f^i_n)\}_{1\leq i\leq
M}$ is a basis of $ V^+\cap z^N U$.
\end{lem}

\begin{pf}
We have to compute the determinant of the inverse image of the complex:
$$\L\,\longrightarrow \, V/V^+$$
by the morphism $\bar \mu_U^N$, which is $\prod
(1-T_j\frac{t_{i j}}z)\big)^{-1}\cdot U\to V/V^+$. It is
straightforward that its determinant coincides with that of the following
complex:
$$\im(h) \,\longrightarrow \, V/z^N U$$
where $h$ is the homothety:
$$\begin{aligned}
h: V^+ &\longrightarrow V^+ \\
(f_1,\ldots, f_n) & \mapsto \big(\prod_{i=1}^N(z-t_{i1})f_1,\ldots
\prod_{i=1}^N(z-t_{in})f_n\big)
\end{aligned}$$

Recall that $V^+=E[[z]]=k[[z]]\oplus\overset{n}\ldots\oplus k[[z]]$.
Consider now the following evaluation map:
{\small $$\begin{aligned}
k[[z]]\oplus\overset{n}\ldots\oplus k[[z]] &\to
{\mathcal M}_{Nn}\,:=\,\bigoplus_{\Sb 1\leq i\leq N\\ 1\leq j\leq
n\endSb}k[[t_{ij}]]
\\ (f_1(z),\ldots, f_n(z))&\mapsto
\big(
 f_1(t_{11}),\ldots,f_1(t_{N1}),\ldots,f_n(t_{1n}),\ldots,f_n(t_{Nn})
\big)
\end{aligned}
$$}
and observe that the cokernel of $h$ is isomorphic to:
$$k[[z]]/\prod_i(z-t_{i1})\oplus\ldots\oplus
k[[z]]/\prod_i(z-t_{in})$$

Now, we construct the following exact sequence of
complexes (written vertically):
{\small $$\minCDarrowwidth19pt\CD
0 @>>> \im(h) @>>> V^+ @>>> \operatorname{Coker}(h) @>>> 0 \\
@. @V{\pi}VV @V{(\pi,v^N)}VV @V{v^N}VV \\
0 @>>> V/z^N U @>>> (V/z^N U)\oplus {\mathcal M}_{Nn}
@>>> {\mathcal M}_{Nn} @>>> 0
\endCD$$}

Since the determinant of the morphism in the complex of the left hand
side is precisely $(\bar \mu^N_U)^*\Omega_+$, it is sufficient to calculate
the others. First, observe that:
$$\operatorname{det}(\bar v_N)=\prod_{j=1}^n\Delta_j\qquad \text{where
}\Delta_j:=\prod_{i<k}(t_{ij}-t_{kj})$$

To compute the determinant of the complex in the middle, note that in the
diagramm:
{\small $$\minCDarrowwidth16pt\CD
0 @>>> V^+\cap z^N U @>>> V^+ @>>> V^+/(V^+\cap z^N U) @>>> 0 \\
@. @VVV @VVV @VVV \\
0 @>>> {\mathcal M}_{Nn} @>>> {\mathcal M}_{Nn}\oplus
V/z^N U @>>> V/z^N U @>>> 0
\endCD$$}
the determinant of the morphism of the right hand side is a non-zero
constant for $N>>0$ (since $U$ lies in $\gr^0(V)$ and $\Omega_+(U)\neq
0$). We conclude now since in the  determinant of the morphism of the
complex of the left hand side is:
$$\operatorname{det}
\pmatrix f^1_1(t_{11}) & \ldots & f^1_1(t_{N1})  & \ldots & f^1_n(t_{1n})
& \ldots & f^1_n(t_{Nn}) \\
\vdots & & & & & & \vdots \\
  f^M_1(t_{11}) & \ldots & f^M_1(t_{N1})  & \ldots & f^M_n(t_{1n})
& \ldots & f^M_n(t_{Nn}) 
\endpmatrix$$
 where $M:=N\cdot n$ and $\{f^i=(f^i_1,\ldots, f^i_n)\vert 1\leq i\leq
M\}$ is a basis of $ V^+\cap z^N U$.
\end{pf}

\begin{rem}
Since the previous Lemma shows that $(\prod_j\Delta_j)\big(
(\mu_U^N)^*\Omega_+\big)(U)$ is a determinant for all $U$, it follows that
one can generalize the ``Addition Formula'' of \cite{SS} to this case.
\end{rem}

\begin{thm}\label{thm:BAmain}
Let $U_0\in\gr^m(V)$ be a rational point. Then,  there exists
polynomials $p_{ij}\in k[\{t_{ij}\}]$ and $\zeta_{m}\in V^+_{m}$
generating $V^+_{m}/V^+_{m-1}$ such that for every rational point, $U$,
lying on a (Zarisky) open neighbourhood of $U_0$ the following formula
holds:
$$\psi_U(z,s)=z\cdot\zeta_m^{-1}\cdot
\sum_{i>0}\big(f_1^i(z)p_{i1}(s),\ldots, f_n^i(z)p_{in}(s)\big)$$
where $\{f^i(z)=(f^i_1(z),\ldots, f^i_n(z))\}_{i>0}$ is a basis of $U$.
\end{thm}

\begin{rem}
Lemma \ref{lem:main} and Theorem \ref{thm:BAmain} have straightforward
generalizations for the case $\dim\T\leq\dim E$.
\end{rem}

\begin{pf}
Once we have proved the Lemma \ref{lem:main}, this proof is a
generalization of that of Theorem 4.8 of \cite{MP}.

Let us consider a rational point $U\in\gr^m(V)$. In order to compute the
$j$-th entry of its BA function, we will calculate first the expression
$\frac{\phi_{j,N}^*\Omega_+}{\phi_N^*\Omega_+}(U)$ for all $N>>0$, where
$\phi_N:\prod^N \w C^n\times\{U\}\to \grv$ and  $\phi_{j,N}$ consist of
composing with the Abel morphism in the $j$-th entry. However, we shall
use the fact that $\phi_{j,N}$ and the morphism $\prod^N \w
C^n\times\{U_{j}\}\to \grv$ coincide ($U_j:= (1,\ldots, (1-\frac{\bar
z}z)^{-1},\ldots,1)$ in the $j$ place).

Observe that $U_{j}\in\gr^{m+1}(V)$ and choose $\zeta_{m}$ such that
$\zeta_m^{-1}\cdot U_0\in\gr^0(V)$ lies on the complementary of the zero
locus of the section $\Omega_+$. Then, it is clear that the points
$U\in\gr^m(V)$ such that $\Omega_+(\zeta_m^{-1}\cdot U)\neq 0$ define a
(Zariski open) neighborhood of $U_0$, and that Lemma \ref{lem:main} might
be applied to $\zeta_m^{-1}\cdot U$.

Let $\{f^1,\ldots,f^{M+1}\}$ be a basis of $\big((1,\ldots,
z^{-1},\ldots,1)z^{-N}V^+\big)\cap \zeta_m^{-1} U$ such that
$\{f^1,\ldots,f^{M}\}$ is a basis of
$z^{-N}V^+\cap \zeta_m^{-1} U$. Let $\bar f^i$ be $z^N\cdot f^i$. Applying
the Lemma \ref{lem:main}, one gets:
$$\frac{\phi_{j,N}^*\Omega_+}{\phi_N^*\Omega_+}(U)\,=\,
\prod_{l=1}^N(\bar z-t_{lj})^{-1}\cdot\bar z\cdot
\big(\bar f^{M+1}_j\prod_l t_{lj}+\sum_{l=1}^M\bar
f^l_j\cdot \bar p_{lj}(t)\big)$$ (up to a scalar). Taking
inverse limit in $N$, replacing $\bar z$ by $z$, and recalling that the
$j$-th entry of the BA function is
$\prod(1-\frac{t_{lj}}z)^{-1}\cdot
\frac{\phi_{j,N}^*\Omega_+}{\phi_N^*\Omega_+}(U)$, it follows that
the $j$-th entry is:
$$z\cdot\sum_{l>0} f^l_j(z)p_{lj}(t)$$ where
$\{f^l=(f^l_1,\ldots,f^l_n)\vert l>0\}$ is a basis of $\zeta_m^{-1} U$.
Finally, the very construction of the polynomials
$p_{lj}(t)$ implies that they are symmetric in the $t$, so that they
can be  expressed  in terms of their symmetric functions $s$. The claim is
proved.
\end{pf}

\subsection{Bilinear Identity}

The previous Theorem together with the definition of the pairing $\res$
imply easily the following:

\begin{thm}\label{thm:bilid}
Let $U,U'\in\grv$ be two rational points, then the Residue Bilinear
Identity holds:
$$\res\big(\frac{1}z\psi_U(z,s)\big)\cdot\big(\frac1z\psi_{U'}^*(z,s')
\big)\,=\,0$$ if and only if $U\subseteq U'$ (the equality holds precisely
when $U$ and $U'$ are in the same connected component).
\end{thm}

\section{$n$-component KP hierarchy}

Here, we introduce the $n$-component KP hierarchy ($n$-KP) in a very
concise way. We will define it as a system of Lax equations and follow
closely \cite{Sa} (see also \cite{vL}). Another common approach is based on
representation theory (see, for instance, \cite{DJKM,KvL}). The last
approach might be ``included'' in ours by studing the action of the linear
group of $E((z))$ on the space of global sections of the determinant
bundle (see \cite{P2}).

\subsection{Pseudodifferential Operators}

Let us begin this section summarizing some standard definitions and
properties of pseudo-differential operators (pdo).

For a $\C$-algebra $A$ and a $\C$-derivation $\partial:A\to A$, one 
considers the $A$-module of pdo: 
$$\P\,:=\,\left\{\sum_{i\leq n}a_i\partial^i\,\vert\, a_i\in A, 
n\in\Z\right\}$$ 
 
The following generalization of the Leibnitz rule: 
$$\big( \sum_i a_i\partial^i\big)\big( \sum_j b_j \partial^j\big) \,:=\, 
\sum_{i,j}\sum_{k\geq 0}\binom{i}{k}a_i(\partial^k b_j)\partial^{i+j-k}$$ 
endows $\P$ with a $\C$-algebra structure. Moreover, $\P$ contains a 
distinguished $\C$-algebra; namely, the algebra $\D$ of differential 
operators (those elements $\sum_{i\leq n}a_i\partial^i$ such that 
$a_i=0$ for all 
$i\leq 0$). 
 
A pdo $\sum_{i\leq n}a_i\partial^i$ is called of 
order $n$ iff $a_n\neq 0$. The subspace of the operators of order less 
or equal than $n\in\Z$ will be denoted by $\P(n)$. Since  $\P=\D\oplus 
\P(-1)$, every operator $P$ decomposes as a sum $P_++P_-$.  Finally, 
define the adjoint of $P=\sum_{i\leq n}a_i\partial^i$ to be 
$P^*=\sum_{i\leq n}(-\partial)^ia_i$. 
 
Observe that: 
\begin{itemize} 
\item $\P(n)\P(m)\subseteq \P(n+m)$;
\item the Leibnitz rule induces a composition law in the affine 
subspace $1+\P(-1)\subset \P$;
\item $1+\P(-1)$ acts transitively on $\partial+\P(-1)$ by conjugation;
\item the stabilizer of $\partial$ consists of those pdo with constant coefficients. (Here, $a$ is constant iff $\partial a=0$).
\end{itemize} 

In the particular case of matrix valued functions, we impose one
constraint; namely, the leading coefficient, $a_n$, is the identity matrix,
$\id$.

\begin{defn}
A $n\times n$-matrix-valued oscillating function is a formal expression of
the type:
$$\left(\id+\sum_{i<0} M_i(s)\bar z^i\right)\cdot
e^{\xi(s,\bar z)}$$
 where $M_i(s)$ are $n\times
n$-matrices, and:
 $$\xi(s,\bar z)\,:=\, \sum_{i>0}\pmatrix s_{i1} & & \\ & \ddots
& \\ & & s_{in}\endpmatrix \bar z^{i}$$
\end{defn}

From now on, we will consider oscillating functions and pdo over the ring $M_{n\times n}(C[[s]])$ ($s=\{s_{ij}\}_{\Sb i>0
\\ 1\leq j\leq n\endSb}$). Define $\partial_{ij}:=\frac{d}{d s_{ij}}$ and  $\partial=\sum_{j=1}^n\partial_{1j}$. 

\subsection{$n$-component KP}

Let $L, C^{1)},\ldots, C^{n)}$ be  pdo with
$n\times n$ matrix coefficients of the form:
$$\begin{aligned}
L\,&=\, \id\partial+L_1(s)\partial^{-1}+
L_2(s)\partial^{-2}+\ldots  \\
C^{(i)}\,&=\,
E_{i}+C^{(i)}_1(s)\partial^{-1}+C^{(i)}_2(s)\partial^{-2}\ldots
\end{aligned}$$
where $E_{i}$ is a matrix whose only non-zero entry is $1$ in the $(i,i)$
place.

Then, the $n$-component KP hierarchy is the following set of Lax
equations:
\beq
\cases
\partial_{ij}L\,=\,[(L^iC^{(j)})_+,L] &\\
\partial_{ij}C^{(k)}\,=\, [(L^i C^{(j)})_+,C^{(k)}]&
\endcases
\qquad 1\leq j,k \leq n\, ,\, i>0
\label{eq:nKPLax}\end{equation}
where $LC^{(j)}=C^{(j)}L$,
$C^{(j)}C^{(k)}=\delta_{jk}C^{(j)}$ and 
$\sum_{j=1}^n C^{(j)}=\id$.

The above system might be regarded as the compatibility condition of the
following system of differential
equations:
\beq
L w=\bar z\cdot w\quad,\quad
C^{(j)}w= w\cdot E_j \quad,\quad
\partial_{ij}w= B^{(j)}_i w
\label{eq:nKPwave}\end{equation}
for a formal oscillating matrix function
$w(\bar z,s)=(\id+\sum_{l<0} w_l(s)\bar z^l)e^{\xi(s,\bar z)}$ and
$B^{(j)}_i:=((L C^{(j)})^i)_+=(L^i C^{(j)})_+$.

\begin{thm}
There are 1-1 correspondence between the set of solutions of the
$n$-component KP hierarchy and the set of rational points
of  the open subscheme of
$\gr^0(V)$ given by $\Omega_+\neq 0$ ($V=M_{n\times n}(k((z)))$ and
$z=\frac1{\bar z}$). 
\end{thm}

For proving the Theorem we shall need the following generalization of the
Lemma of \cite{DJKM} proved in \cite{KvL}:

\begin{lem}
Let $P,Q$ be two pdo with matrix
coefficients. If:
$$\res_{z=0} P(s,\partial)e^{\xi(s,z)}\cdot Q^t(s',\partial')e^{\xi(s',z)} dz
\,=\,0$$
then $(P\cdot Q^*)_-=0$. (Here, the superscript $t$ denotes the
transpose).
\end{lem}

\begin{pf}
It is well known that the system \ref{eq:nKPLax} has a solution $w$ if and
only if there exists a pdo $P=\id+\sum_{i<0}P_i(s)\partial^i$
satisfying \ref{eq:nKPwave} with $w$ replaced by $P$. In that case, it
follows that $w_i(s)=P_i(s)$.

Further, if there exists $P$ as above solving \ref{eq:nKPwave}, then
$L=P\partial P^{-1}$, $C^{(i)}=PE_i P^{-1}$.

Let $P=\id+\sum_{i<0}P_i(s)\partial^i$ be a pdo. If it is a solution, we
have:
$$\partial_{ij}P(s,\partial)\,=\, 
\bar z^iP(s,\partial)E_j-(L^iC^{(j)})_-P(s,\partial)$$
thus, for $w(s,\bar z)=P(s,\bar z)e^{\xi(s,\bar z)}$:
$$\partial_{ij}w(s,\bar z)\,=\, \bar
z^iP(s,\bar z)E_je^{\xi(s,\bar z)}-(L^iC^{(j)})_-w(s,\bar z)$$
and:
$$\partial_{ij}w(s,\bar z)\vert_{s=0}\,=\, \bar
z^i (0,\ldots,\overset{(j)}1,\ldots,0)+\sum_{k<i}\bar z^k w_k$$
where $w_k$ are certain vector valued functions. Summing up, given a
solution $P$ we have shown that the vector space, $U$, generated by
$w(s,\bar z)$ (as the variables $s$ vary) belongs to $\gr^0(V)$. It is
easy to check that $\Omega_+(U)\neq 0$.

Conversely, given a point $U\in \gr^0(V)$ let $P$ be a pdo such that
$\psi_U(z,s)=P(s,\partial)e^{\xi(s,\bar z)}$. Theorems \ref{thm:BAmain} and
\ref{thm:bilid} imply that its BA functions satisfy the following relation:
$$\res\big(\frac1z
(\partial_{ij}-B_i^{(j)})
\psi_U(z,s)\big)\cdot\big(\frac1z\psi_U^*(z,s')\big)
\,=\,0\qquad \forall\, 1\leq j\leq n, i\geq 1
$$
so, the Lemma implies $((\partial_{ij}P-B_i^{(j)}P)\cdot P^*)_-=0$.

Observe that $(\partial_{ij}-B_i^{(j)})P$ is of negative order since:
{\small $$
\begin{aligned}
& (\partial_{ij}-B_i^{(j)})\psi_U(z,s)=
(\partial_{ij}P+z^iPE_j-B_i^{(j)}P)e^{\xi(s,\bar z)}=
\\ & \qquad  =
(\partial_{ij}P+P\partial^iE_j-B_i^{(j)}P)e^{\xi(s,\bar z)}=
(\partial_{ij}P+L^iC^{(j)}P-B_i^{(j)}P)e^{\xi(s,\bar z)}=
\\ & \qquad  =
(\partial_{ij}P+(L^iC^{(j)})_-P)e^{\xi(s,\bar z)}
\end{aligned}$$}

Then, it follows that $\partial_{ij}P-B_i^{(j)}P$ must be identically
zero; that is, $\psi_U$ is a solution of the $n$-KP.
\end{pf}

\begin{rem}
Observe that the equations \ref{eq:nKPLax} determine $P$ up to right
multiplication by  an operator $\id+\sum_{i<0}A_i\partial^i$ with constant
coefficients. (\cite{vL}).
\end{rem}

Let us finish this section with a brief comment on Wronskian solutions.
This ``Wronskian Method'' permits us to construct solutions of the $n$-KP
starting with $n$ solutions of the KP; equivalently (in terms of
Grassmannians), it is a procedure to construct a point of
$\gr(k((z))^{\oplus n^2})$ from a point of $\gr(k((z))^{\oplus n})$.

Consider the ``Wronskian embedding'':
$$ \begin{aligned}{\mathcal W}: \gr(k((z))^{\oplus n})&\hookrightarrow
\gr(k((z))^{\oplus n^2})
\\
U\,&\mapsto U\oplus U^{1)}\oplus\ldots\oplus U^{n-1)}
\end{aligned}$$
where:
$$U^{i)}\,:=\, \{(\frac{d}{dz})^i f(z)\,\vert\, f(z)\in U\}$$
Now, an easy calculation shows the following relation between the BA
functions:
$$\psi_{{\mathcal W}(U)}\,=\,{\operatorname{Wronskian}}(\psi_U)$$
that is, if $\psi_U$ is the vector valued function
$(\psi_U^{(1)},\ldots,\psi_U^{(n)})$ then $\psi_{{\mathcal W}(U)}$ is the
determinant of the matrix $\{(\frac{d}{dz})^{i-1} \psi_U^{(j)}\}_{1\leq
i,j\leq n}$.

\section{Moduli Spaces}

Now, the results of \cite{MP} may be generalized in order to
give equations for the moduli space of algebraic data:
$$\w
\M_{g,n}^r\,=\,\left\{(C,p_1,\ldots,p_n,\alpha_1,\ldots,\alpha_n,M,\beta)\right\}$$
where $C$ is a curve, $p_i\in C$, $\alpha_i:\w\o_{C,p_i}\simeq k[[z]]$,
$M$ is a rank $r$ torsion free sheaf, and $\beta:\w M_{\{p_i\}}\iso
k[[z]]^{\oplus n\cdot r}$. (Here, $\w{\,}$ denotes the completion of a sheaf along a divisor). But, let us be more precise.

Given a flat curve $\pi:C\to S$ and a Cartier divisor $D$ denote:
$$\widehat\o_{C,D} \,=\, \limpl{n}\o_C/\o_C(-n)$$
where $\o_C(-1)$ is the ideal sheaf of $D$. Assume that $D$ is of finite degree, flat and smooth over $S$. (Here smooth over $S$ means that for
every closed point $x\in D$ there exists an open neighborhood $U$ of
$x$ in $C$ such that the morphism $U\to S$ is smooth). Then, $\widehat\o_{C,D}$ is a sheaf of $\o_S$-algebras. We also define the following sheaf of
$\o_S$-algebras:
$$\widehat \Sigma_{C,D}\,=\,\limil{m}\widehat\o_{C,D}(m)$$

\begin{defn}
Let $S$ be a $k$-scheme. Define the functor $\tilde\M_{g,n}^r$
over the category of $k$-schemes by:
$$S\rightsquigarrow \tilde\M_{g,n}^r(S)=
\{\text{ families $(C,p_1,\ldots,p_n,\alpha_1,\ldots,\alpha_n,M,\beta)$
over $S$ }\}$$  where these families satisfy:
\begin{enumerate}
\item $\pi:C\to S$ is a proper flat morphism,
whose geometric fibres are reduced curves of arithmetic genus $g$;
\item $p_i:S\to C$ ($1\leq i\leq n$) is a section of $\pi$, such that 
(when considered as a Cartier Divisor over $C$, also denoted by
$p_i$) is of relative degree 1, flat and smooth over $S$. 
\item for each irreducible component of $C$ there is at least one divisor
$p_i$ lying on it;
\item $\alpha_i$ ($1\leq i\leq n$) is an isomorphism of
$\o_S$-algebras
$\widehat
\Sigma_{C,p_i}\,\iso\, \o_S((z))$.
\item $M$ is a torsion free rank $r$ bundle on $C$.
\item $\beta$ is a direct sum of $n$ isomorphisms of $\o_S$-modules
$\w M_{p_i}\iso\o_S[[z]]^{\oplus r}$ ($1\leq i\leq n$).
\end{enumerate}
\end{defn}

On the set $\tilde\M_{g,n}^r(S)$ one can define an
equivalence  relation, $\sim$: $(C,\{p_i,\alpha_i\},M,\beta)$ and $(C',\{p'_i,\alpha_i'\},M',\beta')$ are said
to be equivalent, if there exists an isomorphism $C\to C'$ (over $S$)
such that the first family goes to the second under the induced
morphisms.

\begin{defn}
The moduli functor of $\widehat\M_{g,n}^r$, is the
functor over the category of
$k$-schemes defined by the sheafication of the functor:
$$S\rightsquigarrow {\tilde\M_{g,n}^r(S)}/{\sim}$$
\end{defn}

\begin{thm}
There is an injective map of functors:
{\small $$\begin{aligned}
\w\M_{g,n}^r &\longrightarrow \gr(k((z))^{\oplus
n})\times\gr(k((z))^{\oplus n\cdot r})
\\ (C,\{p_i,\alpha_i\},M,\beta) &\mapsto
\big(H^0(C-\{p_1,\ldots,p_n\},\o_C),H^0(C-\{p_1,\ldots,p_n\},M)\big)
\end{aligned}$$}

Moreover, $\widehat\M_{g,n}^r$ is representable by a closed subscheme of
$\gr^{1-g}(k((z))^{\oplus
n})\times\gr(k((z))^{\oplus n\cdot r})$.
\end{thm}

\begin{pf}
Let $(C,\{p_i,\alpha_i\},M,\beta)$ be a point of $\widehat\M_{g,n}^r$.
Define the divisor $D$ by $p_1+\ldots+p_n$ and the sheaf $M(m)$ by
$M\otimes\o_C(mD)$. Proceeding as in Proposition~6.3 of \cite{MP} one
shows that the sheaf
$\limil{m}\pi_*M(m)$ is an $S$-valued point of
$\gr(\widehat \Sigma_{C,\{p_i\},}^{\oplus n\cdot r})$.
Moreover, it lies on the component of index $\deg(M)+r(1-g)$.

Now, the morphism of the statement maps $(C,\{p_i,\alpha_i\},M,\beta)\in
\widehat\M_{g,n}^r(S)$ to the pair:
$$\big(\limil{m}\pi_*\o_C(m), \limil{m}\pi_*M(m)\big)$$
which are understood as submodules of $\o_S((z))^{\oplus n}$ and of
$\o_S((z))^{\oplus n\cdot r}$ via $\alpha_1\oplus\ldots\alpha_n$ and
$\beta$, respectively.

It is clear that the image lies on the subset of  $\gr^{1-g}(k((z))^{\oplus
n})\times\gr(k((z))^{\oplus n\cdot r})$ defined by:
\beq
\{(A,B) \text{ such that }A\cdot A=A\text{ and }A\cdot B=B\}
\label{eq:subscheme}\end{equation}
which is actually a closed subscheme (by the same arguments of the proof
of Theorem~6.5 of \cite{MP}). Here, the composition  laws are given as
follows: an element $(f_1,\ldots,f_n)$ of $k((z))^{\oplus n}$ is
represented as the diagonal matrix ${\small \pmatrix f_1 & & 0\\ &
\ddots &
\\ 0 & & f_n\endpmatrix}$; an element of $k((z))^{\oplus n\cdot r}$ is
thought as a $n\times r$ matrix; then, the composition laws are those
induced by the multiplication of matrices).

Starting with a pair $(A,B)$ of that subscheme, a well known procedure
(see, for instance, Theorem~6.4 in \cite{MP}) enables us to construct
algebro-geometric data $(C,\{p_i,\alpha_i\},M,\beta)\in\widehat\M_{g,n}^r$
whose image is $(A,B)$. Further, both constructions are inverse of each
other.

Finally, it is worth observing that our construction, although analogous,
does not follows from the  $n=1$ case, since $H^0(C-\{p_1,\ldots,p_n\})$
is a subalgebra of $k((z))^{\oplus
n}$, while $H^0(C-\{p_1\})\oplus\ldots\oplus
H^0(C-\{p_n\})$ does not.
\end{pf}

\begin{rem}
A direct consequence of the Theorem \ref{thm:BAmain} is that (with a
suitable normalization) the restriction of our BA function
 to a rational point corresponding to
$(C,\{p_i,\alpha_i\},M,\beta)$ gives  the standard notions of Ahkiezer
functions in all its flavours: vector valued, matrix valued,
multi--punctured (\cite{Kr1,Kr3,KrN}, see also
\cite{Kr4,Pr,PW,SW}).

In this case, the Bilinear Residue Identity is geometrically meaningful;
it is equivalent to the fact that the sum of the residues must vanish.
\end{rem}

\begin{thm}\label{thm:equ-mod}
Let $(A,B)$ be a rational point of $\gr^{1-g}(k((z))^{\oplus
n})\times\gr(k((z))^{\oplus n\cdot r})$ and let $\zeta_A\in V$ such that
$\Omega_+(\frac{\zeta_A}{z}\cdot A)\neq 0$. Then $(A,B)$ lies in
$\w\M_{g,n}^r$ if and only if their BA functions satisfy:
$$\cases
\res\big((1,\ldots,1)\big)\cdot
\big(\frac{\zeta_A}{z}\psi^*_A(z,s)\big)\,=\,0 & \\
\res\big(\frac{\zeta_A}{z}\psi_A(z,s)\cdot
\frac1z\psi_A(z,s')\big)\cdot
\big(\frac1z\psi_A^*(z,s'')\big)\,=\,0 & \\
\res\big(\frac{\zeta_A}{z}\psi_A(z,s)\cdot\frac1z\psi_B(z,s')\big)\cdot
\big(\frac1z\psi_B^*(z,s'')\big)\,=\,0 & 
\endcases$$
\end{thm}

\begin{pf}
Firstly, recall that the pairing considered in $k^{\oplus n\cdot r}$
is $(A,B)\,\mapsto \tr(A\cdot B^t)$. Using Theorems \ref{thm:bilid} and
\ref{thm:BAmain}, it is easy to check that the equations defining the
subscheme \ref{eq:subscheme} are those given in the statement.
\end{pf}

Note that the equations given in this Theorem can be transformed into
a set of differential equations for tau functions that might be used to
characterize certain sections of bundles over
$\w\M_{g,n}^r$. (See \cite{MP} for how this transformation is made).
On the other hand, these equations can also be written as algebraic
relations among sections of the determinant bundle.

\section{Solutions to the $n$-KP and non-abelian theta functions}


Fix two positive integers $r,d$ such that $(r,d)=1$. Let us fix a smooth
curve $C$ and a vector bundle $F$ such that $\rank(F)=r$ and 
$\deg(F)=-d+r(g-1)$. 

Denote by $\U_s(r,d)$ the moduli
space of rank $r$ degree $d$ stable vector bundles on
$C$. Recall that $\U_s(r,d)$ carries a
natural line bundle associated to the Weil divisor given by:
$$\{M\text{ such that } h^0(M\otimes F)>0\}$$

The most important properties of this bundle follows from its construction
as a determinant. Let us recall the contruction following \cite{DN,Po}. 

Let $\P$ be a Poincar\'e bundle on $C\times\U_s(r,d)$ and $p_i$ ($i=1,2$)
the projection onto the $i$-th factor. Let us consider a resolution of
$\Q:=\P\otimes p_1^*F$ by locally free sheaves:
$$\CD 0 @>>> V_1 @>\alpha>> V_0 @>>>0 \\
@. @. @VVV \\
& & & &\Q \endCD$$
such that $p_{2*}(V_0)=0$. Since the vertical arrow is a quasi-isomorphism,
it induces an isomorphism:
$$\det(R p_{2*}\Q)\,\simeq\, \det(R p_{2*}V_{\cdot})$$
Note that $\det(R p_{2*}V_{\cdot})\simeq \wedge(R^1p_{2*}(V_0))\otimes
\wedge(R^1p_{2*}(V_1))^{-1}$ (where $\wedge$ denotes the top exterior
algebra). It is now easy to check that the morphism $\alpha$ induces a
section:
$$\theta_F:=det(\alpha)\in H^0(\U_s(r,d), \det(R p_{2*}\Q)^{-1})$$

Lemma~2.1 of \cite{Po} proves that: a) the above construction does not
depend on the resolution; b) the line bundle $\det(R p_{2*}\Q)$ does depend
only on the equivalence class of $F$, $[F]$, in the Grothendieck group of
algebraic coherent sheaves over $C$, $K(C)$ (see also \cite{DN}); and, c)
the zero locus of the section $det(\alpha)$ is the subscheme whose closed
points are those bundles $M$ such that $h^0(C,M\otimes F)\neq 0$. 
Moreover, Lemma~1.2 of \cite{Po} shows that $\det(R p_{2*}\Q)$ does
not depend on the choice of a universal bundle (up to isomorphisms)
provided that $d\rank(F)+r\chi(F)=0$; or, equivalently, $\chi(M\otimes F)=0$
for all $M\in \U_s(r,d)$ (this is why we have taken
$\deg(F)=-d+r(g-1)$).

Summing up, given a class $[F]\in K(C)$ (such that
$d\rank(F)+r\chi(F)=0$) there is an associated line bundle on $\U_s(r,d)$,
$\o(\Theta_{[F]}):=\det(R p_{2*}\Q))$. For each element $F'\in [F]$, there is
a section of the dual of this bundle, $\theta_{F'}$, whose zero locus is $\{M\text{
s.t. }h^0(C,M\otimes F')\neq 0\}$.

Let us now relate the above picture with infinite Grassmannians. As
before, some data must be fixed.

Let $r,d$ be  two positive integers such that
$(r,d)=1$. Let us fix $(C,p,\alpha,F,\gamma)\in\M_{g,1}^{r}$ such that $C$
is smooth and $\deg(F)=-d+r(g-1)$. 

From \S4 we know that:
$$
\w\U(r,d)\,:=\, \{(M,\beta)\text{ s.t. }
(C,p,\alpha,M,\beta)\in\w\M_{g,1}^r\text{ and }\deg(M)=d\}
$$
is a scheme. Let $\w\U_s(r,d)$ be the open subscheme consisting of the
points of $\w\U(r,d)$ such that $M$ is stable. 

If $\M$ is the universal bundle of $\w\U_s(r,d)$ and $V=M_{r\times
r}(k((z)))$, then the submodule:
$$\limil{m} \pi_* (\M(m)\otimes F)\,\hookrightarrow\, V\w\otimes
\o_{\w\U(r,d)}$$
(via $\beta\otimes\gamma$ and $\alpha$) corresponds to the
 morphism of schemes:
$$\iota_{F}:\w\U_s(r,d)   \hookrightarrow \gr^0(V)$$
which for rational points is given by:
$$(M,\beta)   \mapsto H^0(C-p,M\otimes F)$$ 

The above discussion and the exactness of the following
sequence:
$$0\to \pi_* (M\otimes F) \to \limil{m} \pi_* (M(m)\otimes
F) \to V/V^+ \to R^1 \pi_* (M\otimes F) \to 0$$
show the following:

\begin{thm}\label{thm:tau-theta}
Let $\pi:\w\U_s(r,d)\to\U_s(r,d)$ be the canonical projection. It holds
that:
\begin{itemize}
\item $(\iota_{F}^*\Omega_+)(M,\beta)\neq 0\,\iff\,
h^0(C-p,M\otimes F)\neq0$;
\item $\iota_{F}^*\Omega_+=\pi^*\theta_F$ (up to a non zero constant);
\item $\iota_{F}^*\det\iso\pi^*\o(\Theta_{[F]})$.
\end{itemize}
\end{thm}

Roughly said, the Theorem states that, analogously the the rank 1 case
(see \cite{Kr2,Sh}), it is theoretically possible to give solutions for
the $n$-KP in terms of certain non-abelian theta functions. 
However, further
research must be made in this direction to obtain explicit expressions.




\end{document}